\documentclass[11pt]{article}
\usepackage{graphicx}
\usepackage{epsfig}
\usepackage{array}
\usepackage{amssymb}
\usepackage{picinpar}
\usepackage{hhline}
\usepackage{amsmath} 
\usepackage{mathrsfs} 
\usepackage{wrapfig}
\usepackage{indentfirst}
\usepackage[numbers]{natbib}
\usepackage{enumitem}
\usepackage{float}
\def\rsq{\hspace*{\fill}$\Box$}
\newtheorem{prop}{Proposition}[section]

\newtheorem{thm}[prop]{Theorem}

\newtheorem{cor}[prop]{Corollary}

\newtheorem{ex}{Example}[section]

\def\pf{\noindent {\bf Proof: }}

\def\D{\Delta}

\def\s{\sigma}
\def\l{\lambda}

\def\di{\displaystyle}

\def\g{\gamma}

\def\Ra{\Rightarrow}

\makeindex             

\title{Some Results on incidence coloring, star arboricity and domination number}

\author{Pak Kiu Sun\thanks{The research were partially supported by the Pentecostal Holiness
Church Incorporation (Hong Kong)}\\
\small Department of Mathematics\\[-0.8ex]
\small Hong Kong Baptist University\\[-0.8ex]
\small Kowloon Tong, Hong Kong\\
\small \texttt{lionel@hkbu.edu.hk}\\
\and
Wai Chee Shiu\\
\small Department of Mathematics\\[-0.8ex]
\small Hong Kong Baptist University\\[-0.8ex]
\small Kowloon Tong, Hong Kong\\
\small \texttt{wcshiu@hkbu.edu.hk} }

\date{}

\begin{document}
\maketitle
\begin{abstract}
Two inequalities bridging the three isolated graph invariants,
incidence chromatic number, star arboricity and domination number,
were established. Consequently, we deduced an upper bound and a
lower bound of the incidence chromatic number for all graphs. Using
these bounds, we further reduced the upper bound of the incidence
chromatic number of planar graphs and showed that cubic graphs with
orders not divisible by four are not 4-incidence colorable. The
incidence chromatic numbers of Cartesian product, join and union of
graphs were also determined.

\end{abstract}


\section{Introduction}\label{sec-int}
 An incidence coloring separates the whole graph into disjoint
independent incidence sets. Since incidence coloring was introduced
by Brualdi and Massey \cite{BR93}, most of the researches were
concentrated on determining the minimum number of independent
incidence sets, also known as the incidence chromatic number, which
can cover the graph. The upper bound of the incidence chromatic
number of planar graphs \citep{DO04}, cubic graphs \citep{MA05} and
a lot of other classes of graphs were determined \citep{DO04, DO05,
SH02, SUN08, WA02}. However, for general graphs, the best possible
upper bound is an asymptotic one \citep{GU97}. Therefore, to find an
alternative upper bound and lower bound of the incidence chromatic
number for all graphs is the main objective of this paper.

In Section~\ref{sec-star}, we will establish a global upper bound
for the incidence chromatic number in terms of chromatic index and
star arboricity. This result reduces the upper bound of the
incidence chromatic number of the planar graphs. Also, a global
lower bound which involves the domination number will be
introduced in Section~\ref{sec-dom}. Finally, the incidence
chromatic number of graphs constructed from smaller graphs will be
determined in Section~\ref{sec-product}.

All graphs in this paper are connected. Let $V(G)$ and $E(G)$ (or
$V$ and $E$) be the vertex-set and edge-set of a graph $G$,
respectively. Let the set of all neighbors of a vertex $u$ be
$N_G(u)$(or simply $N(u)$). Moreover, the degree $d_G(u)$(or
simply $d(u)$) of $u$ is equal to $|N_G(u)|$ and the maximum
degree of $G$ is denoted by $\D(G)$ (or simply $\D$).  All
notations not defined in this paper can be found in the books
\cite{BO76} and \cite{WE01}.

Let $D(G)$ be a digraph induced from $G$ by splitting each edge
$e(u,v) \in E(G)$ into two opposite arcs $uv$ and $vu$. According
to \cite{SH02}, incidence coloring of $G$ is equivalent to the
coloring of arcs of $D(G)$. Two distinct arcs $uv$ and $xy$ are
{\it adjacent} provided
one of the following holds:\medskip\\
(1) $u=x$; \\
(2) $v=x$ or $y=u$.\medskip

Let $A(G)$ be the set of all arcs of $D(G)$. An {\it incidence
coloring} of $G$ is a mapping $\s : A(G) \to C$, where $C$ is a
{\it color-set}, such that adjacent arcs of $D(G)$ are assigned
distinct colors. The {\it incidence chromatic number}, denoted by
$\chi_i$, is the minimum cardinality of $C$ for which $\s : A(G)
\to C$ is an {\it incidence coloring}. An {\it independent set} of
arcs is a subset of $A(G)$ which consists of non-adjacent arcs.

\section{Incidence chromatic number and Star
arboricity}\label{sec-star}

A {\it star forest} is a forest whose connected components are
stars. The {\it star arboricity} of a graph $G$ (introduced by
Akiyama and Kano \cite{AK85}), denoted by $st(G)$, is the smallest
number of star forests whose union covers all edges of $G$.

We now establish a connection among the chromatic index, the star
arborcity and the incidence chromatic number of a graph. This
relation, together with the results by Hakimi et al. \cite{HA96},
provided a new upper bound of the incidence chromatic number of
planar graphs, $k$-degenerate graphs and bipartite graphs.

\begin{thm}\label{thm-upper} Let $G$ be a graph. Then $\chi_i(G) \leq \chi'(G) +
st(G)$, where $\chi'(G)$ is the chromatic index of $G$. \end{thm}

\pf We color all the arcs going into the center of a star by the
same color. Thus, half of the arcs of a star forest can be colored
by one color. Since $st(G)$ is the smallest number of star forests
whose union covers all edges of $G$, half of the arcs of $G$ can
be colored by $st(G)$ colors. The uncolored arcs now form a
digraph which is an orientation of $G$. We color these arcs
according to the edge coloring of $G$ and this is a proper
incidence coloring because edge coloring is more restrictive.
Hence $\chi'(G) + st(G)$ colors are sufficient to color all the
incidences of $G$. \rsq

\medskip

We now obtain the following new upper bounds of the incidence
chromatic numbers of planar graphs, a class of $k$-degenerate
graphs and a class of bipartite graphs.

\begin{cor}\label{cor-plane} Let $G$ be a planar graph. Then $\chi_i(G) \leq \D +
5$ for $\D \neq 6$ and $\chi_i(G) \leq 12$ for $\D = 6$.
\end{cor} \pf The bound is true for $\D \leq
 5$, since Brualdi and Massey \cite{BR93} proved that $\chi_i(G) \leq 2\D$. Let $G$ be a planar graph with $\D
 \geq 7$, we have $\chi'(G) = \D$ \citep{SA01,VI68}. Also,
Hakimi et al. \cite{HA96} proved that $st(G) \leq 5$. Therefore,
$\chi_i(G) \leq \D + 5$ by Theorem~\ref{thm-upper}. \rsq

While we reduce the upper bound of the incidence chromatic number of
planar graphs from $\D + 7$ \citep{DO04} to $\D +5$, Hosseini Dolama
and Sopena \cite{DO05} reduced the bound to $\D + 2$ under the
additional assumptions that $\D \geq 5$ and girth $g \geq 6$.

A {\it $k$-degenerate graph} $G$ is a graph with vertices ordered
$v_1, v_2, \dots, v_n$ such that each $v_i$ has degree at most $k$
in the graph $G[v_1,v_2,\dots,v_i]$. A {\it restricted
$k$-degenerate graph} is a $k$-degenerate graph with the graph
induced by $N(v_i) \cap \{v_1,v_2,\dots,v_{i-1}\}$ is complete for
every $i$. It has been proved by Hosseini Dolama et al.\cite{DO04}
that $\chi_i(G) \leq \D + 2k -1$, where $G$ is a $k$-degenerate
graph. We lowered the bound for restricted $k$-degenerate graph as
follow.

\begin{cor}\label{cor-kdeg} Let $G$ be a restricted
$k$-degenerate graph. Then $\chi_i(G) \leq \D + k +2$.
\end{cor} \pf By Vizing's theorem, we have $\chi'(G) \leq \D +1$.
Also, the star arboricity of a restricted $k$-degenerate graph $G$
is less than or equal to $k+ 1$ \citep{HA96}. Hence we have
$\chi_i(G) \leq \D + k +2$ by Theorem~\ref{thm-upper}.

\begin{cor}\label{cor-bip} Let $B$ be a bipartite graph with at most one cycle.
Then $\chi_i(B) \leq \D + 2$. \end{cor} \pf Hakimi et al.
\cite{HA96} proved that $st(B) \leq 2$ where $B$ is a bipartite
graph with at most one cycle.  Also, it is well known that $\chi'(B)
= \D$. These results together with Theorem~\ref{thm-upper} proved
the corollary. \rsq

\section{Incidence chromatic number and Domination number}\label{sec-dom}
A {\it dominating set} $S \subseteq V(G)$ of a graph $G$ is a set
where every vertex not in $S$ has a neighbor in $S$. The {\it
domination number}, denoted by $\g(G)$, is the minimum cardinality
of a domination set in $G$.

A {\it maximal star forest}  is a star forest with maximum number
of edges. Let $G=(V,E)$ be a graph, the number of edges of a
maximal star forest of $G$ is equal to $|V| - \g(G)$ \citep{FE02}.
We now use the domination number to form a lower bound of the
incidence chromatic number of a graph. The following proposition
reformulates the ideas in \cite{AL89} and \cite{MA05}.

\begin{prop}\label{prop-dom} Let $G = (V,E)$ be a graph. Then $\di \chi_i(G) \geq
\frac{2|E|}{|V| - \g(G)}$.\end{prop} \pf Each edge of $G$ is
divided into two arcs in opposite directions. The total number of
arcs of $D(G)$ is therefore equal to $2|E|$. According to the
definition of the adjacency of arcs, an independent set of arcs is
a star forest. Thus, a maximal independent set of arcs is a
maximal star forest. As a result, the number of color class
required is at least $\di \frac{2|E|}{|V| - \g(G)}$. \rsq

\begin{cor}\label{cor-dom} Let $G = (V,E)$ be an $r$-regular graph.
Then $\di \chi_i(G) \geq \frac{r}{1 - \frac{\g(G)}{|V|}}$.
\end{cor} \pf By Handshaking lemma, we have $\di 2|E| = \sum_{v\in V}
d(v) = r|V|$, the result follows from Proposition~\ref{prop-dom}.
\rsq

Corollary~\ref{cor-dom} provides an alternative method to show
that a cycle $C_n$, where $n$ is not divisible by 3, is not
3-incidence colorable. As $C_n$ is a 2-regular graph with $\di
\g(C_n) > \frac{|V(C_n)|}{3}$,  we have  $$\di \chi_i(C_n) \geq
\frac{2}{1-\frac{g(G)}{|V|}} > \frac{2}{1-\frac{1}{3}} = 3.$$

\begin{cor}\label{cor-nec} Let $G = (V,E)$ be an $r$-regular graph.
Two necessary conditions for  $\di \chi_i(G) = r +1$ (also for
$\chi(G^2) = r +1$ \cite{SUN08}) are:\begin{enumerate} \item The
number of vertices of $G$ is divisible by $r+1$. \item If $r$ is
odd, then the chromatic index of $G$ is equal to $r$.
\end{enumerate}
\end{cor} \pf We prove 1 only, 2 was proved in \cite{SUN06}. By
Corollary~\ref{cor-dom}, if $G$ is  an $r$-regular graph and $\di
\chi_i(G) = r+1$, then $\di r+1 = \chi_i(G) \geq \frac{r|V|}{|V|-
\g(G)} \Ra \frac{|V|}{r+1} \geq \g(G)$. Since the global lower
bound of domination number is $\di \left\lceil \frac{|V|}{\D +
1}\right\rceil$, we conclude that the number of vertices of $G$
must be divisible by $r +1$. \rsq
%

\section{Graphs Constructed from Smaller Graphs}\label{sec-product}

In this section, we determine the upper bound of the incidence
chromatic number of union of graphs, Cartesian product of graphs
and join of graphs, respectively. Also, these bounds can be
attained by some classes of graphs \citep{SUN07}. Let the set of
colors assigned to the arcs going into $u$ be $C^+_G(u)$.
Similarly, $C^-_G(u)$ represents the set of colors assigned to the
arcs going out from $u$.

We start by proving the following theorem about union of graphs.

%

\begin{thm}\label{thm-union} For all graphs $G_1$ and $G_2$, we have
$\chi_i(G_1 \cup G_2) \leq \chi_i(G_1)+ \chi_i(G_2)$.
\end{thm}
\pf If some edge $e \in E(G_1) \cap E(G_2)$,  then we delete it
from either one of the edge set. This process will not affect
$I(G_1 \cup G_2)$, hence, we assume $E(G_1) \cap E(G_2) =
\varnothing$. Let $\s$ be a $\chi_i(G_1)$-incidence coloring of
$G_1$ and $\l$ be a $\chi_i(G_2)$-incidence coloring of $G_2$
using different color set. Then $\phi$ is a proper $(\chi_i(G_1) +
\chi_i(G_2))$-incidence coloring of $G_1 \cup G_2$ with $\phi(uv)
= \s(uv)$ when $e(uv) \in E(G_1)$ and  $\phi(uv) = \l(uv)$ when
$e(uv) \in E(G_2)$. \rsq
\medskip

The following example revealed that the upper bound given in
Theorem~\ref{thm-union} is sharp.

\begin{ex} {\rm Let $n$ be an even integer and not divisible by 3.
Let $G_1$ be a graph with $V(G_1) = \{u_1, u_2, \dots, u_n\}$ and
$E(G_1) = \{u_{2i-1}u_{2i}|1\leq i \leq \frac{n}{2}\}$.
Furthermore, let $G_2$ be another graph with $V(G_2) = V(G_1)$ and
$E(G_2) = \{u_{2i}u_{2i+1}| 1\leq i \leq \frac{n}{2}\}$. Then, it
is obvious that $\chi_i(G_1) = \chi_i(G_2) = 2$ and $G_1 \cup G_2
= C_n$ where $n$ is not divisible by 3. Therefore, $\chi_i(C_n) =
4 =\chi_i(G_1) + \chi_i(G_2)$.} \rsq \end{ex}

Next, we prove the theorem about the Cartesian product of graphs.
The following definition should be given in prior.

\begin{thm}\label{thm-carte} For all graphs $G_1$ and $G_2$, we have
$\chi_i(G_1 \Box G_2) \leq \chi_i(G_1)+ \chi_i(G_2)$.
\end{thm}
\pf Let $|V(G_1)| =m$ and $|V(G_2)| =n$. $G_1 \Box G_2$ is a graph
with $mn$ vertices and two types of edges: from conditions (1) and
(2) respectively. The edges of type (1) form a graph consisting of
$n$ disjoint copies of $G_1$, hence  its incidence chromatic
number equal to $\chi_i(G_1)$. Likewise, the edges of type (2)
form a graph with incidence chromatic number  $\chi_i(G_2)$.
Consequently, the graph $G_1 \Box G_2$ is equal to the union of
the graphs from (1) and (2). By Theorem~\ref{thm-union}, we have
$\chi_i(G_1 \Box G_2) \leq \chi_i(G_1)+ \chi_i(G_2)$. \rsq

We demonstrate the upper bound given in Theorem~\ref{thm-carte} is
sharp by the following example.

\begin{ex} {\rm Let $G_1 = G_2 = C_3$, then $G_1 \Box G_2$ is a 4-regular graph. If it
is 5-incidence colorable, then its square has chromatic number
equal to 5 \citep{SUN08}. However, all vertices in $G_1 \Box G_2$
is of distance at most 2. Therefore, $G_1 \Box G_2$ is not
5-incidence colorable and the bound derived in
Theorem~\ref{thm-carte} is attained.
 } \rsq \end{ex}

Finally, we consider the incidence chromatic number of the join of
graphs.

\begin{thm}\label{thm-join} For all graphs $G_1$ and $G_2$ with $|V(G_1)| = m
$, $|V(G_2)| = n $ and $m \geq n \geq 2$. We have $\chi_i(G_1 \vee
G_2) \leq \min\{m + n, \max\{\chi_i(G_1), \chi_i(G_2)\} + m +
2\}$.
\end{thm}
\pf On the one hand, we have $\chi_i(G_1 \vee G_2) \leq m + n$. On
the other hand, the disjoint graphs $G_1$ and $G_2$ can be colored
by $\max\{\chi_i(G_1),\chi_i(G_2)\}$ colors, and all other arcs in
between can be colored by $m+2$ new colors. Therefore,
$\max\{\chi_i(G_1), \chi_i(G_2)\} + m + 2$ is also an upper bound
for $\chi_i(G_1 \vee G_2)$. \rsq

Similar to the previous practices , we utilize the following
example to show that the upper bound in Theorem~\ref{thm-join} is
sharp.

\begin{ex} {\rm Let $G_1 = K_m$ and $G_2 = K_n$. Then the upper bound $m+n$
is obtained since   $G_1 \vee G_2 \cong K_{m+n}$. On the other
hand, let  $G_1$ be the null graph of order $m$ and  $G_2$ be the
null graph of order $n$. Then the other upper bound
$\max\{\chi_i(G_1), \chi_i(G_2)\} + m + 2$ is attained because
$G_1 \vee G_2 \cong K_{m,n}$.} \rsq
\end{ex}




\begin{thebibliography}{10}

\bibitem{AK85}
Jin Akiyama and Mikio Kano, \emph{Path factors of a graph}, Graphs
and
  applications : proceedings of the First Colorado Symposium on Graph Theory
  (Frank Harary and John~S. Maybee, eds.), Wiley, 1985, pp.~1--21.

\bibitem{AL89}
I.~Algor and N.~Alon, \emph{The star arboricity of graphs}, Discrete
Math.
  \textbf{75} (1989), 11--22.

\bibitem{BO76}
J.~A. Bondy and U.~S.~R. Murty, \emph{Graph theory with
applications}, 1st ed.,
  New York: Macmillan Ltd. Press, 1976.

\bibitem{BR93}
R.~A. Brualdi and J.~J.~Q. Massey, \emph{Incidence and strong edge
colorings of
  graphs}, Discrete Math. \textbf{122} (1993), 51--58.

\bibitem{FE02}
Sheila Ferneyhough, Ruth Haas, Denis Hanson, and Gary MacGillivray,
\emph{Star
  forests, dominating sets and \textsc{R}amsey-type problems.}, Discrete Math.
  \textbf{245} (2002), 255--262.

\bibitem{GU97}
B.~Guiduli, \emph{On incidence coloring and star arboricity of
graphs},
  Discrete Math. \textbf{163} (1997), 275--278.

\bibitem{HA96}
S.~L. Hakimi, J.~Mitchem, and E~Schmeichel, \emph{Star arboricity of
graphs},
  Discrete Math. \textbf{149} (1996), 93--98.

\bibitem{DO05}
M.~Hosseini~Dolama and E.~Sopena, \emph{On the maximum average
degree and the
  incidence chromatic number of a graph}, Discrete Math. and Theoret. Comput.
  Sci. \textbf{7} (2005), 203--216.

\bibitem{DO04}
M.~Hosseini~Dolama, E.~Sopena, and X.~Zhu, \emph{Incidence coloring
of
  $k$-degenerated graphs}, Discrete Math. \textbf{283} (2004), 121--128.

\bibitem{MA05}
M.~Maydanskiy, \emph{The incidence coloring conjecture for graphs of
maximum
  degree 3}, Discrete Math. \textbf{292} (2005), 131--141.

\bibitem{SA01}
Daniel~P. Sanders and Yue Zhao, \emph{Planar graphs of maximum
degree seven are
  class 1}, J. Combin. Theory Ser. B \textbf{83} (2001), 201--212.

\bibitem{SH02}
W.~C. Shiu, P.~C.~B. Lam, and D.~L. Chen, \emph{Note on incidence
coloring for
  some cubic graphs}, Discrete Math. \textbf{252} (2002), 259--266.

\bibitem{SUN06}
Wai~Chee Shiu, Peter Chi~Bor Lam, and Pak~Kiu Sun, \emph{Cubic
graphs with
  different incidence chromatic numbers}, Congr. Numerantium \textbf{182}
  (2006), 33--40.

\bibitem{SUN08}
Wai~Chee Shiu and Pak~Kiu Sun, \emph{Invalid proofs on incidence
coloring},
  Discrete Math. \textbf{308} (2008), 6575--6580.

\bibitem{SUN07}
Pak~Kiu Sun, \emph{Incidence coloring: Origins, developments and
relation with
  other colorings}, Ph.D. thesis, Hong Kong Baptist University, 2007.

\bibitem{VI68}
V.~G. Vizing, \emph{Some unsolved problems in graph theory (in
{R}ussian)},
  Uspekhi Mat Nauk \textbf{23} (1968), 117--134.

\bibitem{WA02}
S.~D. Wang, D.~L. Chen, and S.~C. Pang, \emph{The incidence coloring
number of
  \textsc{H}alin graphs and outerplanar graphs}, Discrete Math. \textbf{256}
  (2002), 397--405.

\bibitem{WE01}
D.~B. West, \emph{Introduction to graph theory}, 2nd ed.,
Prentice-Hall, Inc.,
  2001.

\end{thebibliography}

\end{document}